\documentclass{gtart}
\usepackage{graphicx}
\usepackage{indentfirst}
\usepackage{amsmath,amsfonts,amsthm,amssymb}
\usepackage{mathrsfs}
\usepackage{amscd}

\newcommand{\homo}{{\mathrm H}}
\newcommand{\relspin}{\underline{\mathrm{Spin}}^c}
\newcommand{\dif}{{\mathrm d}}
\def\co{\colon\thinspace}

\newtheorem{thm}{Theorem}[section]
\newtheorem{lem}[thm]{Lemma}
\newtheorem{cor}[thm]{Corollary}
\newtheorem{prop}[thm]{Proposition}

\newtheorem*{thm*}{Theorem}

\theoremstyle{definition}
\newtheorem{defn}[thm]{Definition}
\newtheorem{rem}[thm]{Remark}

\begin{document}

\title{Link Floer homology detects the Thurston norm}
\authors{Yi Ni}
\address{Department of Mathematics, Princeton University, Princeton,  New Jersey 08544}

\begin{abstract}
We prove that, for a link $L$ in a rational homology 3--sphere,
the link Floer homology detects the Thurston norm of its
complement. This generalizes the previous results due to
Ozsv\'ath, Szab\'o and the author.
\end{abstract}


\maketitlepage

\section{Introduction}

Link Floer homology was introduced by Ozsv\'ath and Szab\'o
\cite{OSz7}, as a multi-filtered theory for links in rational
homology 3--spheres. This theory generalizes an earlier invariant
for knots, the knot Floer homology \cite{OSz4} \cite{Ra}.

One interesting topic in Floer theory is the relationship with the
Thurston norm. For knot (and link) Floer homology, this topic was
studied for links in integer homology 3--spheres in \cite{OSz5},
\cite{Ni} and \cite{OSz8}. In particular, Ozsv\'ath and Szab\'o
showed that, for a link $L\subset S^3$, the link Floer homology
detects the Thurston norm of the complement of $L$. Although not
stated explicitly, their proof actually works for links in integer
homology spheres.

In the current paper, we will generalize Ozsv\'ath and Szab\'o's
result to links in rational homology 3--spheres.

In Subsection \ref{sectfil}, we will define an affine function
$$\mathfrak H\co\relspin(Y,L)\to\homo^2(Y,L;\mathbb Q).$$ Then the link Floer homology provides a function
$$y\co\homo_2(Y,L;\mathbb Q)\to\mathbb Q,$$
defined by
$$y(h)=\max_{\{\mathfrak r\in\relspin(Y,L)|\:\widehat{HFL}(Y,L,\mathfrak r)\ne0\}}\langle \mathfrak H(\mathfrak r),h\rangle.$$

\begin{thm}\label{LinkNorm}
Suppose $L$ is an oriented link in a rational homology 3--sphere
$Y$, $M=Y-\mathrm{int}(Nd(L))$. Suppose $h\in\homo_2(M,\partial
M;\mathbb Q)$ is an integral class, then $h$ can be represented by
a properly embedded surface without sphere components. Let
$\chi(h)$ be the maximal possible value of the Euler
characteristic of such surfaces. Then
$$-\chi(h)+\sum_{i=1}^l|h\cdot[\mu_i]|=2y(h).$$
Here $\mu_1,\dots,\mu_l$ are the meridians of the components of
$L$.
\end{thm}

\begin{rem}\label{Thurston}
The term $-\chi(h)$ is almost the Thurston norm of $h$ \cite{T}.
In fact, if the boundary tori of $M$ are all incompressible, then
we can rewrite the equality in the above theorem as
$$x(h)+\sum_{i=1}^l|h\cdot[\mu_i]|=2y(h),$$
where here $x(\cdot)$ is the Thurston norm.
\end{rem}

\begin{rem}
Suppose $M$ is a compact 3--manifold with boundary consisting of
tori, and $\homo_2(M)=0$. Then $M$ is the complement of a link in
a rational homology sphere. Theorem \ref{LinkNorm} gives a
criterion to determine whether any component of $\partial M$ is
compressible, in the terms of link Floer homology. If $T$ is a
torus in a rational homology 3--sphere $Y$, then $T$ splits $Y$
into two rational homology solid tori. Thus Theorem \ref{LinkNorm}
 also gives a criterion to determine whether $T$ is
compressible.

Incompressible tori play a very important role in ``traditional
3--dimensional topology". We hope that the above observation will
be useful for studying the relationship between Floer homology and
traditional 3--dimensional topology.
\end{rem}

The paper is organized as follows. Section 2 contains a rather
general result about the existence of longitudinal foliations.
This result will be the starting point of our proof. Then, in
Section 3, we generalize the main result in \cite{Ni} to
null-homologous oriented links in rational homology spheres. In
Section 4, we give some preliminaries on link Floer homology. In
particular, we discuss the relative Spin$^c$ structures. Section 5
will be devoted to the proof of the main theorem. We use the
``cabling trick" from \cite{OSz8}, as well as the techniques from
\cite{He}, to reduce the general case of our main theorem to the
case proved in Section 3.

\noindent{\bf Acknowledgements.} We are grateful to David Gabai
and Zolt\'an Szab\'o for some helpful conversations. We also wish
to thank Efstratia Kalfagianni for informing Kaiser's work
\cite{Kai}.

The author is partially supported by the Centennial fellowship of
the Graduate School at Princeton University.

\section{Longitudinal foliations}

As in \cite{OSz5} and \cite{Ni}, when one tries to relate Floer
homology with Thurston norm, the first step is always to establish
an existence result about taut foliations. In this section, we are
going to establish the corresponding result we need.

\begin{defn}
$M$ is an $n$--dimensional manifold. A smooth (codimension--$1$)
{\it foliation} $\mathscr F$ of $M$ is a smooth integrable
hyperplane field on $M$. A {\it leaf} $L$ of $\mathscr F$ is a
maximal path-connected integral submanifold for $\mathscr F$.
\end{defn}

By abuse of notation, we also denote the collection of the leaves
by $\mathscr F$.

From now on, we assume the foliation is {\it co-oriented}, namely,
there exists a unit vector field on the manifold, which is
transverse to the foliation everywhere.

If $\gamma$ is a path in a leaf $L$, then $\mathscr F$ defines a
parallel transport in a small neighborhood of $\gamma$. Let $a,b$
be the two ends of $\gamma$, there are two small transversals
$I_a,I_b$ passing through $a,b$, so that the parallel transport
along $\gamma$ defined by $\mathscr F$ gives a diffeomorphism of
$I_a$ onto $I_b$. Moreover, if $\gamma$ is a loop with base point
$b$, then the germ of the diffeomorphism at $b$ is called the {\it
holonomy} along the loop $\gamma$.

Every closed orientable 3--manifold admits a smooth foliation. So
in order to extract useful topological information about
3--manifolds out of foliations, one needs some further restriction
on the foliations.

\begin{defn}
Let $\mathscr F$ be a foliation of a 3--manifold $M$. $\mathscr F$
is {\it taut} if there exists a closed curve intersecting every
leaf of $\mathscr F$ transversely.
\end{defn}

In order to study knot Floer homology, one always needs some
additional conditions on the taut foliations. For example, the
foliations should be ``longitudinal". The definition is as
follows.

\begin{defn}
Suppose $K\subset Y$ is a null-homologous knot, $\mathscr F$ is a
taut foliation of $Y-\mathrm{int}(Nd(K))$. We say that $\mathscr
F$ is
 a {\it longitudinal foliation}, if the restriction of  $\mathscr F$ on $\partial
Nd(K)$ consists of longitudes.
\end{defn}

Gabai shows that longitudinal foliations exist in many cases,
including the classical knots \cite{G1}. In fact, we are going to
prove the following rather general result about the existence of
longitudinal foliations.

\begin{prop}\label{LongiFol}
Suppose $K\subset Y$ is a null-homologous knot, $Y-K$ is
irreducible. Then for any fibred knot $J\subset S^3$ with
sufficiently large genus, $Y-\mathrm{int}(Nd(K\#J))$ admits a
smooth longitudinal foliation with a compact leaf, which is the
minimal genus Seifert surface of $K\#J$.
\end{prop}

Proposition \ref{LongiFol} is a weak form of the following theorem
due to Gabai.

\begin{thm*}[Gabai, \cite{G2}]
Suppose $K_i\subset Y_i$ are nontrivial null-homologous knots,
$Y_i-K_i$ are irreducible, $i=1,2$. Then
$Y_1\#Y_2-\mathrm{int}(Nd(K_1\#K_2))$ admits a smooth longitudinal
foliation with a compact leaf, which is the minimal genus Seifert
surface of $K_1\#K_2$.\hfill\qedsymbol
\end{thm*}

\begin{proof}[Proof of Proposition \ref{LongiFol}]
Suppose $F\subset M=Y-\mathrm{int}(Nd(K))$ is a minimal genus
Seifert surface for $K$. By the main theorem in \cite{G}, there
exists a taut smooth foliation $\mathscr F$ of $M$, so that

(1) $\mathscr F\pitchfork\partial M$, and $\mathscr F|\partial M$
has no Reeb component,

(2) $F$ is a leaf of $\mathscr F$,

(3) if $\theta$ is a closed curve in $F$,
$f\co(-\varepsilon,\varepsilon)\to(-\delta,\delta)$ is a
representative of the germ of the holonomy along $\theta$, then
$$
\frac{\dif^kf}{\dif t^k}(0)=\left\{\begin{array}{lr}
1,&k=1,\\0,&k>1.
\end{array}\right.
$$

Here the above (3) holds by the Induction Hypothesis (iii) in the
proof of \cite[Theorem 5.1]{G}.

Cut $M$ open along $F$, we get a sutured manifold
$(M_0,\gamma_0)$, and $\mathscr F$ becomes a foliation $\mathscr
F_0$ of $M_0$. The suture $\gamma_0$ is an annulus. By the above
condition (1), $\mathscr F_0|\gamma_0$ is determined by a global
holonomy $f\co I\to I$. Namely, pick the square $I\times I$,
foliated by $I\times t$'s. Glue $0\times I$ with $1\times I$ by a
diffeomorphism $f$, then the induced foliation on $S^1\times I$ is
equivalent to the foliation $\mathscr F_0|\gamma_0$. We can view
$\gamma_0$ as the union of two squares $a\times I$ and $b\times
I$, so that the restriction of the foliation in $a\times I$
consists of $a\times t$'s, and the holonomy takes place in
$b\times I$.

Suppose $D_8$ is an octagon with edges
$a_1,b_1,a_2,b_2,\dots,a_4,b_4$ in cyclic order. Consider
$D_8\times I$, foliated by $D_8\times t$'s. Let $g,h\co I\to I$ be
two diffeomorphisms with the two ends fixed. We glue $b_1\times I$
with $b_3\times I$ by the map $\mathrm{id}\times g$, glue
$b_2\times I$ with $b_4\times I$ by the map $\mathrm{id}\times h$.
The new manifold is $R\times I$, with an induced foliation
$\mathscr G$. Here $R$ is a genus 1 compact surface with one
boundary component. Obviously, $\mathscr G|\partial R\times I$ has
a global holonomy $[g,h]$. We can view $\partial R\times I$ as the
union of two squares $a'\times I$ and $b'\times I$, so that the
restriction of the foliation in $a'\times I$ consists of $a'\times
t$'s, and the holonomy takes place in $b'\times I$.

Now we glue the two sutured manifolds $(M_0,\gamma_0)$ and
$(R\times I,\partial R\times I)$ together, so that $a\times I$ is
glued to $a'\times I$ by the identity. Then the new sutured
manifold $(M_1,\gamma_1)$ has an induced foliation $\mathscr F_1$,
so that $\mathscr F_1|\gamma_1$ has a global holonomy
$f\circ[g,h]$.

Repeat the above construction $m$ times, we get a foliated sutured
manifold $(M_m,\gamma_m)$, which is the union of $(M_0,\gamma_0)$
and $(R_m\times I,\partial R_m\times I)$ along a square in the
suture, and the holonomy of the foliation on $\gamma_m$ is
$$f\circ[g_1,h_1]\circ[g_2,h_2]\circ\cdots\circ[g_m,h_m].$$
Here $R_m$ is a compact genus $m$ surface with one boundary
component. Now we can make use of the following theorem.

\begin{thm*}[Mather--Sergeraert--Thurston, \cite{Se}, see also
\cite{G}] If $f\co I\to I$ is a $C^{\infty}$ map satisfying
$$
\frac{\dif^kf}{\dif t^k}(\alpha)=\left\{\begin{array}{lr}
1,&k=1,\\0,&k>1,
\end{array}\right.
$$
for $\alpha\in\{0,1\}$, then there exist $C^{\infty}$
diffeomorphisms $g_i,h_i\co I\to I$,$i=1,\dots,n$, satisfying the
above conditions so that
$$f\circ[g_1,h_1]\circ[g_2,h_2]\circ\cdots\circ[g_n,h_n]={\mathrm i\mathrm
d}.$$
\end{thm*}

Hence when $m\ge n$, one can choose the holonomies
$g_i,h_i,i=1,\dots,m$, so that the holonomy of $\mathscr
F_m|\gamma_m$ is the identity, thus $\mathscr F_m|\gamma_m$
consists of closed curves.

Now suppose $J\subset S^3$ is a fibred knot with genus $m$, $G$ is
a minimal genus Seifert surface of $J$. Consider the knot $K\#J$,
with Seifert surface $F'$, which is the boundary connected sum of
$F$ and $G$. If we cut $Y-\mathrm{int}(Nd(K\#J))$ open along $F'$,
the sutured manifold we get is just $(M_m,\gamma_m)$. So
$Y-\mathrm{int}(Nd(K\#J))$ admits a smooth longitudinal foliation
with $F'$ being a compact leaf.
\end{proof}

\section{Genera of Links in rational homology spheres}

In this section, we are going to follow the procedure in \cite{Ni}
to generalize the main result there to null-homologous links in
rational homology 3--spheres.

\begin{thm}\label{EulerBound} Suppose $L$ is a null-homologous oriented link in a closed 3--manifold $Z$, $\homo_1(Z;\mathbb Q)=0$.
$|L|$ denotes the number of components of $L$, and $\chi(L)$
denotes the maximal Euler characteristic of the Seifert surfaces
bounded by $L$. Then
$$\frac{|L|-\chi(L)}2=\max\{i|\;\widehat{HFK}(Z,L,i)\ne0\}.$$
\end{thm}

Let $L$ be a null-homologous oriented $l$--component link in a
rational homology 3--sphere $Z$. Ozsv\'ath and Szab\'o define a
knot $\kappa(L)\subset \kappa(Z)$, where here $\kappa(Z)$ is
obtained by adding $l-1$ 3--dimensional tubes $R_1,\dots,R_{l-1}$
to $Z$. Suppose $P_i$ is the belt sphere of the tube $R_i$. The
knot $\kappa(L)$ intersects $P_i$ in exactly 2 points, we can
remove two disks from $P_i$ at these two points, then glue in a
long and thin (2--dimensional) tube along an arc in $\kappa(L)$,
so as to get a torus $T_i$. $T_i$ is homologous to $P_i$, but
disjoint from $\kappa(L)$. These tori generate
$\homo_2(\kappa(Z)-\kappa(L);\mathbb Z)$.

Let $(Y,K)=(\kappa(Z),\kappa(L))$. Let $G$ be a minimal genus
Seifert surface of $K$, and $Y_0$ be the manifold obtained by
0--surgery on $K$. By \cite[Remark~3.2]{Ni}, we can assume $G$ is
obtained by adding $l-1$ bands to a Seifert surface $F$ of $L$
with maximal Euler characteristic. Hence $\chi(G)=\chi(F)-(l-1)$.
Let $\widehat G$ be the extension of $G$ in $Y_0$ obtained by
gluing a disk to $G$.

\begin{prop}\label{sumwithfibred}Let $L$ be a null-homologous oriented link
in a rational homology 3--sphere $Z$, with irreducible complement.
After doing connected sum with some fibred knots in $S^3$, we get
a new link $L^*$. We consider $(Y^*,K^*)=(\kappa(Z),\kappa(L^*))$,
and the 0--surgered space $Y_0^*$. The conclusion is: for a
suitably chosen $L^*$, $Y^*_0$ can be embedded into a closed
symplectic 4--manifold $(X,\Omega)$, so that
$X=X_1\cup_{Y^*_0}X_2$, $b_2^+(X_j)>0$, and
$$\int_{T_i^*}\Omega=0$$
for all $i$. Moreover, $$ \langle c_1(\mathfrak
k(\Omega)),[\widehat{G^*}]\rangle=2-2g(\widehat{G^*}).$$\qed
\end{prop}

Having Proposition \ref{LongiFol}, the proof of Proposition
\ref{sumwithfibred} is the same as the proof of \cite[Proposition
3.12]{Ni}. So we just omit it here.

We also state the following lemma without giving the proof, since
its proof is not different from the proof of  \cite[Lemma
4.1]{Ni}.

\begin{lem}\label{0surgery} $(Y,K)$ is as before. Let $d$ be
an integer satisfying $\widehat{HFK}(Y,K,i)=0$ for $i\ge d$, and
suppose that $d>1$. Then
$$HF^+(Y_0,[d-1])=0,$$
where
$$HF^+(Y_0,[d-1])=\bigoplus_{\langle c_1(\mathfrak s),[\widehat G]\rangle=2(d-1)}HF^+(Y_0,\mathfrak
s).$$\qed
\end{lem}

\begin{proof}[Proof of Theorem \ref{EulerBound}](Compare the proof of  \cite[Theorem 1.1]{Ni})\quad Suppose $L_1,L_2$
are null-homologous oriented links in $Z_1,Z_2$, respectively. We
have
\begin{eqnarray*}
\widehat{HFK}(Z_1,L_1)\otimes\widehat{HF}(Z_2)&\cong&\widehat{HFK}(Z_1\#Z_2,L_1),\\
\widehat{HFK}(Z_1\#Z_2,L_1\#L_2)\otimes\widehat{HF}(S^2\times
S^1)&\cong&\widehat{HFK}(Z_1\#Z_2,L_1\sqcup L_2).
\end{eqnarray*}
By the above formulas, we can assume $Z-L$ is irreducible. Now
apply Proposition \ref{sumwithfibred} to get a symplectic
4--manifold $(X,\Omega)$, $X=X_1\cup_{Y_0^*} X_2$, with
$b_2^+(X_j)>0$, $\int_{T_i^*}\Omega=0$, and
$$ \langle
c_1(\mathfrak
k(\Omega)),[\widehat{G^*}]\rangle=\chi(\widehat{G^*})<0.$$

The sum
\begin{equation}\label{SpinSum}
\sum_{\eta\in\homo^1(Y_0^*)}\Phi_{X,\mathfrak
k(\Omega)+\delta\eta}
\end{equation}
 is calculated by a homomorphism
which factors through $HF^+(Y_0^*\:,\mathfrak
k(\Omega)|_{Y_0^*})$.

$\homo^1(Y_0^*)\cong \mathbb Z^l$ is generated by the Poincar\'e
duals of $[T_1^*],[T_2^*],\dots,[T_{l-1}^*]$ and
$[\widehat{G^*}]$. So the Spin$^c$ structures in (\ref{SpinSum})
are precisely
$$\mathfrak
k(\Omega)+\sum_{i=1}^{l-1}a_i\mathrm{PD}([T_i^*])+b\:\mathrm{PD}([\widehat{G^*}])
\quad(a_i,b\in\mathbb Z).$$ Here $\mathrm{PD}$ is the Poincar\'e
duality map in $X$. The first Chern classes of these Spin$^c$
structures are
$$c_1(\mathfrak
k(\Omega))+2\sum_{i=1}^{l-1}a_i\mathrm{PD}([T_i^*])+2b\:\mathrm{PD}([\widehat{G^*}]).$$

By the degree shifting formula, the degrees of the terms in
(\ref{SpinSum}) are
\begin{eqnarray*}
&&\frac{(c_1(\mathfrak
k(\Omega))+2\sum a_i\mathrm{PD}([T_i^*])+2b\:\mathrm{PD}([\widehat{G^*}]))^2-2\chi(X)-3\sigma(X)}4\\
&=&\frac{(c_1(\mathfrak
k(\Omega))^2-2\chi(X)-3\sigma(X)}4+\sum a_i\langle c_1(\mathfrak k(\Omega)),[T_i^*]\rangle+b\langle c_1(\mathfrak k(\Omega)),[\widehat{G^*}]\rangle\\
&=&\frac{(c_1(\mathfrak k(\Omega))^2-2\chi(X)-3\sigma(X)}4 +
b\chi(\widehat{G^*}).
\end{eqnarray*}

Since $\chi(\widehat{G^*})\ne0$, the terms which have the same
degree as $\Phi_{X,\mathfrak k(\Omega)}$ are precisely those
correspond to $\mathfrak k(\Omega)+\sum a_i\mathrm{PD}([T_i^*])$.
By \cite[Theorem 1.1]{OSz3} and the fact that
$\int_{T_i^*}\Omega=0$, $\Phi_{X,\mathfrak k(\Omega)}$ is the only
nontrivial term at this degree. So $HF^+(Y_0^*\:,\mathfrak
k(\Omega)|_{Y_0^*})$ is nontrivial. Now apply Lemma
\ref{0surgery}, we get our desired result for $L^*$.

The result for $L$ holds by the connected sum
formula.
\end{proof}

As a corollary, we have

\begin{cor}
Suppose $Z$ is a rational homology 3--sphere, $L_+,L_-,L_0\subset
Z$ are 3 null-homologous oriented links, which differ at a
crossing as in the skein relation. Then two of the three numbers
$$\chi(L_+),\chi(L_-),\chi(L_0)-1,$$
are equal and not larger than the third.
\end{cor}
\begin{proof}
In the local picture of the skein relation, if the two strands in
$L_-$ belong to the same component, then $|L_0|=|L_+|+1$, and
there is a surgery exact triangle relating $\widehat{HFK}(Z,L_-)$,
$\widehat{HFK}(Z,L_+)$ and $\widehat{HFK}(Z,L_0)$ \cite{OSz4}. If
$\chi(L_+)<\chi(L_-)$, then
$\widehat{HFK}(Z,L_-,\frac{|L_+|-\chi(L_+)}2)=0$, hence
$$\widehat{HFK}(Z,L_+,\frac{|L_+|-\chi(L_+)}2)\cong\widehat{HFK}(Z,L_0,\frac{|L_+|-\chi(L_+)}2).$$
It follows from Theorem \ref{EulerBound} that
$$|L_0|-\chi(L_0)=|L_+|-\chi(L_+),$$
so $\chi(L_+)=\chi(L_0)-1$. Similarly, one can show that if
$\chi(L_-)<\chi(L_+)$, then $\chi(L_-)=\chi(L_0)-1$, and if
$\chi(L_-)=\chi(L_+)$, then $\chi(L_-)\le\chi(L_0)-1$.

If the two strands in $L_-$ belong to different components, then
$|L_0|=|L_-|-1$, and there is an exact triangle relating
$\widehat{HFK}(Z,L_-)$, $\widehat{HFK}(Z,L_+)$ and
$\widehat{HFK}(Z,L_0)\otimes V$. Here $V=V_{-1}\oplus V_0\oplus
V_{+1}$, $V_{\pm1}\cong\mathbb Z$ are supported in filtration
level $\pm1$, and $V_0\cong\mathbb Z\oplus\mathbb Z$ is supported
in filtration level $0$. An argument similar to the one in the
last paragraph gives the desired result.
\end{proof}

The above result was first proved for links in $S^3$ by
Scharlemann and Thompson \cite{ST}. Then Kaiser proved a much more
general theorem for links in irreducible rational homology
3--spheres \cite{Kai}. Kalfagianni also proved Scharlemann and
Thompson's result for certain links in irreducible homology
3--spheres, and applied it to study the convergence of the HOMFLY
power series link invariants in \cite{Kal}.

\section{Preliminaries on link Floer homology}

In \cite{OSz7}, Ozsv\'ath and Szab\'o defined link Floer homology
for oriented links in rational homology 3--spheres. We will
briefly review the definition and some basic properties.

\subsection{Relative {\rm Spin$^c$} structures}

Let $M$ be a compact 3--manifold with boundary consisting of tori.
There is a canonical isotopy class of translation invariant vector
fields on the torus. Let $v_1$ and $v_2$ be two nowhere vanishing
vector fields on $M$, whose restriction on each component of
$\partial M$ is the canonical translation invariant vector field.
We say $v_1$ and $v_2$ are {\it homologous}, if they are homotopic
in the complement of a ball in $M$. The homology classes of such
vector fields are called {\it relative {\rm Spin}$^c$ structures}
on $M$, and the set of all relative Spin$^c$ structures is denoted
by $\relspin(M,\partial M)$. $\relspin(M,\partial M)$ is an affine
space over $\homo^2(M,\partial M)$.

When $L$ is an oriented link in a closed oriented 3--manifold $Y$,
let $M=Y-\mathrm{int}(Nd(L))$. Then we also denote
$\relspin(M,\partial M)$ by $\relspin(Y,L)$.

There is a natural involution
$$J\co\relspin(M,\partial M)\to\relspin(M,\partial M).$$
If $\mathfrak r\in\relspin(M,\partial M)$ is represented by a
vector field $v$, then $J(\mathfrak r)$ is represented by the
vector field $-v$. Given $\mathfrak r$, $\mathfrak r-J(\mathfrak
r)$ is an element in $\homo^2(M,\partial M)$, denoted by
$c_1(\mathfrak r)$.

\subsection{Heegaard diagrams and {\rm Spin}$^c$
structures}\label{sectdg}

Suppose $L$ is an oriented link in a rational homology sphere
$Y^3$,
$(\Sigma,\mbox{\boldmath${\alpha}$},\mbox{\boldmath${\beta}$},\mathbf
w,\mathbf z)$ is a (generic) balanced $2l$--pointed Heegaard
diagram associated to the pair $(Y,L)$. There is a map
$$\underline{\mathfrak s}_{\mathbf w,\mathbf z}\co\mathbb T_{\alpha}\cap\mathbb T_{\beta}\to\relspin(Y,L),$$
defined in \cite{OSz7}. We sketch the definition of
$\underline{\mathfrak s}_{\mathbf w,\mathbf z}$ as follows.

Let $f\co Y\to[0,3]$ be a Morse function corresponding to the
Heegaard diagram, $\nabla f$ is the gradient vector field
associated to $f$. Let $\gamma_{\mathbf w}$ be the union of the
flowlines of $\nabla f$, such that each of these flowlines
 passes through a point in $\mathbf w$, and connects
an index 0 critical point to an index 3 critical point. Similarly,
define $\gamma_{\mathbf z}$. Suppose $\mathbf x\in\mathbb
T_{\alpha}\cap\mathbb T_{\beta}$, then
$\mbox{\boldmath${\gamma}$}_{\mathbf x}$ denotes the union of the
flowlines connecting index 1 critical points to index 2 critical
points, and passing through the points in $\mathbf x$.

We construct a nowhere vanishing vector field $v$. Outside a
neighborhood $Nd(\gamma_{\mathbf w}\cup\gamma_{\mathbf
z}\cup\mbox{\boldmath${\gamma}$}_{\mathbf x})$, $v$ is identical
with $\nabla f$. Then one can extend $v$ over the balls
$Nd(\mbox{\boldmath${\gamma}$}_{\mathbf x})$. We can also extend
$v$ over $Nd(\gamma_{\mathbf w}\cup\gamma_{\mathbf z})$, so that
the closed orbits of $v$, which pass through points in $\mathbf w$
and $\mathbf z$, give the link $L$. There may be many different
choices to extend $v$ over $Nd(\gamma_{\mathbf
w}\cup\gamma_{\mathbf z})$, we choose the extension as in
\cite[Figure~2]{OSz7}.

Now we let $\underline{\mathfrak s}_{\mathbf w,\mathbf z}(\mathbf
x)$ be the relative Spin$^c$ structure given by $v$. It is easy to
check that $\underline{\mathfrak s}_{\mathbf w,\mathbf z}$ is a
well defined map.

\subsection{Link Floer homology}

Let $\mathbb F_2$ be the field consisting of 2 elements. For
$\mathfrak r\in\relspin(Y,L)$, $\widehat{CFL}(Y,L,\mathfrak r)$ is
a chain complex over $\mathbb F_2$, generated by the $\mathbf x$'s
with $\underline{\mathfrak s}_{\mathbf w,\mathbf z}(\mathbf
x)=\mathfrak r$, and the differential counts holomorphic disks
with $n_{\mathbf w}(\phi)=n_{\mathbf z}(\phi)=0$. The homology of
$\widehat{CFL}(Y,L,\mathfrak r)$ is denoted by
$\widehat{HFL}(Y,L,\mathfrak r)$. And the link Floer homology is
$$\widehat{HFL}(Y,L)=\bigoplus_{\mathfrak r\in\relspin(Y,L)}\widehat{HFL}(Y,L,\mathfrak r).$$

$\widehat{HFL}$ enjoys certain symmetries. In particular, as in
 \cite[Proposition~8.2]{OSz7}, we have

\begin{lem}\label{sJs}
Let $L$ be an oriented link in a rational homology sphere $Y$,
$\mu_1,\dots,\mu_l$ denote the meridians of the components of $L$.
Then
$$\widehat{HFL}(Y,L,\mathfrak
r)\cong\widehat{HFL}(Y,L,J(\mathfrak r )+\sum_{i=1}^l\text{\rm
PD}[\mu_i]).$$\hfill\qedsymbol
\end{lem}

\subsection{An Alexander $\mathbb Q^l$--grading}\label{sectfil}

With the notation as above, we define a function
$$\mathfrak H\co\relspin(Y,L)\to\homo^2(Y,L;\mathbb Q)$$
as follows. Given $\mathfrak r\in\relspin(Y,L)$, let
$$\mathfrak H(\mathfrak r)=\frac{c_1(\mathfrak r)-\sum_{i=1}^l\mathrm{PD}([\mu_i])}2.$$
Moreover, if $\mathbf x\in\mathbb T_{\alpha}\cap\mathbb
T_{\beta}$, we define
$$\mathfrak h_{\mathbf w,\mathbf z}(\mathbf x)=\mathfrak H(\underline{\mathfrak s}_{\mathbf w,\mathbf z}(\mathbf x)).$$

Given $\mathbf x,\mathbf y\in\mathbb T_{\alpha}\cap\mathbb
T_{\beta}$, there exists a closed curve $\omega(\mathbf x,\mathbf
y)\subset \mathbb T_{\alpha}\cup\mathbb T_{\beta}$. $\omega$ is
the union of a curve $a\subset\mathbb T_{\alpha}$ which connects
$\mathbf x$ to $\mathbf y$, and a curve $b\subset\mathbb
T_{\beta}$ which connects $\mathbf y$ to $\mathbf x$. $\omega$ can
also be viewed as a curve in $\Sigma$.

Since $Y$ is a rational homology sphere, there exists a positive
integer $k$, so that $k\omega(\mathbf x,\mathbf y)$ is homologous
to the sum of some copies of $\alpha$ and $\beta$ curves. Let
$\mathcal D$ be such a homology. The following elementary lemma is
important.

\begin{lem}\label{frakH}
With the notation as above, given $\mathbf x,\mathbf y\in\mathbb
T_{\alpha}\cap\mathbb T_{\beta}$, we have
$$\mathfrak h_{\mathbf w,\mathbf z}(\mathbf x)-\mathfrak h_{\mathbf w,\mathbf z}(\mathbf y)
=\frac1k\sum_{i=1}^l(n_{z_i}(\mathcal D)-n_{w_i}(\mathcal
D))\text{\rm PD}[\mu_i].$$
\end{lem}
\begin{proof}
We cap off the copies of $\alpha$ and $\beta$ curves in
$\partial\mathcal D$ to get a 2--dimensional chain $G\subset Y$,
so that $\partial G=k\omega(\mathbf x,\mathbf y)$. $G\cap
(Y-\mathrm{int}(Nd(L)))$ is a homology between $k\omega(\mathbf
x,\mathbf y)$ and some copies of $\mu_i$'s. And the coefficients
of $\mu_i$'s can be computed by counting the algebraic
intersection numbers of $K_i$ with $\mathcal D$. Since
$\underline{\mathfrak s}_{\mathbf w,\mathbf z}(\mathbf
x)-\underline{\mathfrak s}_{\mathbf w,\mathbf z}(\mathbf
y)=\mathrm{PD}([\omega(\mathbf x,\mathbf y)])$ \cite[Lemma
3.11]{OSz7}, we have that
\begin{eqnarray*}
&&k(\mathfrak h_{\mathbf w,\mathbf z}(\mathbf x)-\mathfrak
h_{\mathbf w,\mathbf z}(\mathbf y))\\
&=&\frac k2(c_1(\underline{\mathfrak s}_{\mathbf w,\mathbf
z}(\mathbf x))-c_1(\underline{\mathfrak s}_{\mathbf w,\mathbf
z}(\mathbf y)))\\
&=&k\mathrm{PD}([\omega(\mathbf x,\mathbf y)])\\
&=&\mathrm{PD}(\sum_{i=1}^l(n_{z_i}(\mathcal D)-n_{w_i}(\mathcal
D))[\mu_i]).
\end{eqnarray*}
Hence the result holds.
\end{proof}

The above lemma indicates that $\mathfrak H$ defines a $\mathbb
Q^l$--grading on $\widehat{CFL}(Y,L)$. Following Rasmussen
\cite{Ra}, we call this grading an {\it Alexander grading}. Given
$h^* \in\homo^2(Y,L;\mathbb Q)$, let
$$\widehat{CFL}(Y,L,h^*)\cong\bigoplus_{\mathfrak r\in\relspin(Y,L),\mathfrak H(\mathfrak r)=h^*}\widehat{CFL}(Y,L,\mathfrak r).$$
Then Lemma \ref{sJs} implies that
\begin{equation}\label{h-h}
\widehat{CFL}(Y,L,h^*)\cong\widehat{CFL}(Y,L,-h^*).
\end{equation}

\subsection{A formula for split links}

The following formula for split links will be used in the proof of
Theorem~\ref{LinkNorm}.

\begin{prop}\label{splitlink}
Suppose $L_1\subset Y_1$, $L_2\subset Y_2$ are two oriented links
in rational homology spheres, then $L=L_1\sqcup L_2$ is a split
link in $Y=Y_1\#Y_2$. Let $\mathfrak r_1\in\relspin(Y_1,L_1)$,
$\mathfrak r_2\in\relspin(Y_2,L_2)$, then they naturally give a
relative Spin$^c$ structure $\mathfrak r=\mathfrak r_1\#\mathfrak
r_2 \in\relspin(Y,L)$. We have the following formula:
$$\widehat{CFL}(Y,L,\mathfrak r)\cong\widehat{CFL}(Y_1,L_1,\mathfrak r_1)
\otimes\widehat{CFL}(Y_2,L_2,\mathfrak r_2)\otimes\widehat{HF}(S^1\times S^2).$$
\end{prop}
\begin{proof}
Suppose $L_i$ has $l_i$ components, $i=1,2$. Let
$(\Sigma_i,\mbox{\boldmath${\alpha}$}_i,\mbox{\boldmath${\beta}$}_i,\mathbf
w_i,\mathbf z_i)$ be a weakly admissible balanced $2l_i$--pointed
Heegaard diagram associated to the pair $(Y_i,L_i)$. We construct
a Heegaard diagram for $(Y,L)$ as follows.

Let $A=S^1\times[-1,1]$ be a tube, $\alpha_0=S^1\times0$ is a belt
circle, and $\beta_0$ is a small Hamiltonian perturbation of
$\alpha_0$. Let $\Sigma=\Sigma_1\#\Sigma_2$, with $A$ as the neck
of the connected sum. We put the feet of this tube into two
regions which contain base points. It is easy to verify that
$$(\Sigma,\mbox{\boldmath${\alpha}$}_1\cup\{\alpha_0\}\cup\mbox{\boldmath${\alpha}$}_2,
\mbox{\boldmath${\beta}$}_1\cup\{\beta_0\}\cup\mbox{\boldmath${\beta}$}_2,\mathbf
w_1\cup\mathbf w_2,\mathbf z_1\cup\mathbf z_2)$$ is a weakly
admissible Heegaard diagram for $(Y,L)$.

Now the desired formula can be proved by a standard argument.
\end{proof}

\section{Proof of the main theorem}

In this section, we are going to prove our main theorem. The idea
of the proof is the same as in \cite{OSz8}, but we will take a
slightly different approach.

First of all, let us check Theorem \ref{LinkNorm} for certain
knots in lens spaces. As in \cite{OSz6}, if one does $\frac pq$
surgery on one component of the Hopf link, then the other
component gives a knot $O_{p/q}$ in the lens space $L(p,q)$. The
complement of $O_{p/q}$ is a solid torus, with a meridian disk
$D_{p/q}$. Our result is

\begin{lem}\label{LensKnot}
There are exactly $p$ relative Spin$^c$ structures satisfying that
$$\widehat{HFK}(L(p,q),O_{p/q},\mathfrak r)\ne0.$$
One can denote these relative Spin$^c$ structures by $\mathfrak
r_1,\dots,\mathfrak r_p$, so that
$$\langle c_1(\mathfrak r_i),[D_{p/q}]\rangle=2i-1.$$
Hence Theorem \ref{LinkNorm} holds for $O_{p/q}$.
\end{lem}
\begin{proof}
$(L(p,q),O_{p/q})$ admits a genus 1 Heegaard diagram, such that
$\mathbb T_{\alpha}\cap\mathbb T_{\beta}$ has exactly $p$
intersection points, which correspond to $p$ different relative
Spin$^c$ structures. As in \cite[Lemma 7.1]{OSz6}, we can denote
these relative Spin$^c$ structures by $\mathfrak
r_1,\dots,\mathfrak r_p$, such that $\mathfrak r_{i+1}-\mathfrak
r_i$ is the positive generator of
$\homo^2(L(p,q),O_{p/q})\cong\mathbb Z$, for $i=1,\dots,p-1$. Let
$a_i=\langle c_1(\mathfrak r_i),[D_{p/q}]\rangle$, then
$a_{i+1}-a_i=2$.

Since
$$\langle c_1(\mathfrak r)+c_1(J(\mathfrak r)+\mathrm{PD}[\mu]),[D_{p/q}]\rangle=\langle2\mathrm{PD}[\mu],[D_{p/q}]\rangle=2p,$$
by Lemma \ref{sJs}, the set $\{a_1,a_2,\dots,a_p\}$ admits an
involution $a\mapsto2p-a$. Hence we must have $a_i=2i-1$.

Now it is easy to check Theorem \ref{LinkNorm} for $O_{p/q}$.
\end{proof}

Suppose $L$ is an oriented link in a rational homology 3--sphere
$Y$\hspace{-0.4pt},
\hspace{-0.4pt}$(\Sigma,\mbox{\boldmath${\alpha}$},\mbox{\boldmath${\beta}$},\mathbf
w,\mathbf z)$ is a (generic) balanced $2l$--pointed Heegaard
diagram associated to the pair $(Y,L)$. Given an integral class
$h\in\homo_2(Y,L;\mathbb Q)$, let
$$\mathcal F_{\mathbf w,\mathbf z}^h(\mathbf x)=\langle\mathfrak h_{\mathbf w,\mathbf z}(\mathbf x),h\rangle,$$
for any $\mathbf x\in\mathbb T_{\alpha}\cap\mathbb T_{\beta}$.
Thus $\mathcal F_{\mathbf w,\mathbf z}^h$ defines a $\mathbb
Q$--grading on $\widehat{CFL}(Y,L)$.

\begin{prop}\label{LinkNorm1}
Suppose $L$ is a null-homologous oriented link, and $F$ is a
minimal genus Seifert surface of $L$. Then Theorem \ref{LinkNorm}
holds for $h=[F]$.
\end{prop}
\begin{proof}
As in \cite{OSz7}, we can get a Heegaard diagram
$$(\Sigma',\mbox{\boldmath${\alpha'}$},\mbox{\boldmath${\beta'}$},
w_1,z_l)$$ for $(\kappa(Y),\kappa(L))$, by adding tubes with feet
at $z_i$ and $w_{i+1}$, for $i=1,\dots,l-1$. Suppose $\mathcal D$
is a topological disk in $\mathrm{Sym}^{g+l-1}(\Sigma')$,
$\partial\mathcal D\subset\mathbb T_{\alpha}\cup\mathbb
T_{\beta}$, then $n_{z_i}(\mathcal D)=n_{w_{i+1}}(\mathcal D)$.
Hence
$$n_{z_1}(\mathcal D)-n_{w_l}(\mathcal
D)=\sum_{i=1}^l(n_{z_i}(\mathcal D)-n_{w_i}(\mathcal D)).$$ By
Lemma \ref{frakH}, we conclude that the $\mathbb Q$--grading
defined by $\mathcal F_{\mathbf w,\mathbf z}^{[F]}$ coincides with
the usual Alexander $\mathbb Z$--grading defined on
$\widehat{CFK}(\kappa(Y),\kappa(L))$, as relative gradings.

The proof of  \cite[Theorem 1.1]{OSz7} shows that
$$\widehat{HFL}(Y,L)\cong\widehat{HFK}(Y,L)$$
as relative $\mathbb Q$ graded $\mathbb F_2$--vector spaces.
Moreover, by Lemma~\ref{sJs}, $\widehat{HFL}(Y,L)$, equipped with
the absolute $\mathbb Q$--grading given by $\mathcal F_{\mathbf
w,\mathbf z}^{[F]}$, is symmetric with respect to the origin $0$.
Hence this absolute $\mathbb Q$--grading is identical to the usual
absolute Alexander $\mathbb Z$--grading on $\widehat{HFK}(Y,L)$.

Now we can apply Theorem \ref{EulerBound} to get the conclusion.
\end{proof}

In order to reduce the general case to the case of $h=[F]$, we are
going to use the ``cabling trick" introduced in \cite{OSz8}. The
idea is to consider a $(\mathbf p,\mathbf q)$--cable of $L$. The
method of dealing with cables comes from Hedden's work \cite{He}.

Suppose $L$ is an $l$--component oriented link in $Y$, the
components of $L$ are denoted by $K_1,\dots,K_l$. Let
$(\Sigma,\mbox{\boldmath${\alpha}$},\mbox{\boldmath${\beta}$},\mathbf
w,\mathbf z)$ be a $2l$--pointed Heegaard diagram associated to
the pair $(Y,L)$, satisfying the following conditions:

(1) For each $i\in\{1,2,\dots,l\}$, $\beta_i$ represents a
meridian for $K_i$, namely, $w_i$ and $z_i$ lie on a curve
$\lambda_i$ which meets $\beta_i$ in a single point, and is
disjoint from all the other $\beta$ curves.

(2) $\beta_i$ meets $\alpha_i$ transversely in a single point, and
is disjoint from all the other $\alpha$ curves.

The curve $\lambda_i\subset \Sigma$ is isotopic to the knot $K_i$
in $Y$, hence $\Sigma$ gives a frame on $K_i$.

Suppose $\mathbf p=(p_1,\dots,p_l)$, $\mathbf q=(q_1,\dots,q_l)$
are two $l$--tuples of positive integers, where
$$q_i=p_in_i+1$$
for some $l$--tuple of positive integers $\mathbf
n=(n_1,\dots,n_l)$. We replace $\beta_i$ with a new curve
$\gamma_i$, gotten by performing a ``finger move" of $\beta_i$
along $\lambda_i$ with multiplicity $(p_i-1)$, and then winding
$n_i$ times parallel to $\beta_i$. We put a new basepoint $z_i'$
inside the end of the finger. The new diagram
$(\Sigma,\mbox{\boldmath${\alpha}$},\mbox{\boldmath${\gamma}$},\mathbf
w,\mathbf z')$ gives the link $C(L)=C_{\mathbf p,\mathbf q}(L)$,
which is the $(\mathbf p,\mathbf q)$--cable of $L$ with respect to
the frame specified by $\Sigma$. We can also find a basepoint
$t_i$ outside the finger, so that
$(\Sigma,\mbox{\boldmath${\alpha}$},\mbox{\boldmath${\gamma}$},\mathbf
w,\mathbf t)$ describes $L$. See Figure~1 for an illustration of
the local diagram.

Let $j\co Y-C(L)\to Y-L$ be the inclusion map,
$j_*\co\homo_1(Y-C(L))\to\homo_1(Y-L)$,
$j^*\co\homo^2(Y,L)\to\homo^2(Y,C(L))$ be the induced maps on
(co)homologies. Let $\mu_i'$ be the meridian of
$C_{p_i,q_i}(K_i)$. Then $j_*([\mu_i])=p_i[\mu_i']$.

Consider the intersection points in $\mathbb T_{\alpha}\cap\mathbb
T_{\gamma}$. If the $\gamma_i$--component of a point lies in a
regular neighborhood of $\beta_i$, then we call the intersection
point an {\it $i$--exterior intersection point}. If the
$\gamma_i$--component is supported in a regular neighborhood of
$\lambda_i$, then we call the intersection point an {\it
$i$--interior intersection point}. An intersection point which is
$i$--exterior for all $i=1,\dots,l$ is called an {\it exterior
intersection point}.

\begin{center}
\begin{picture}(375,330)
\put(0,0){\scalebox{0.7}{\includegraphics*[24pt,270pt][559pt,
740pt]{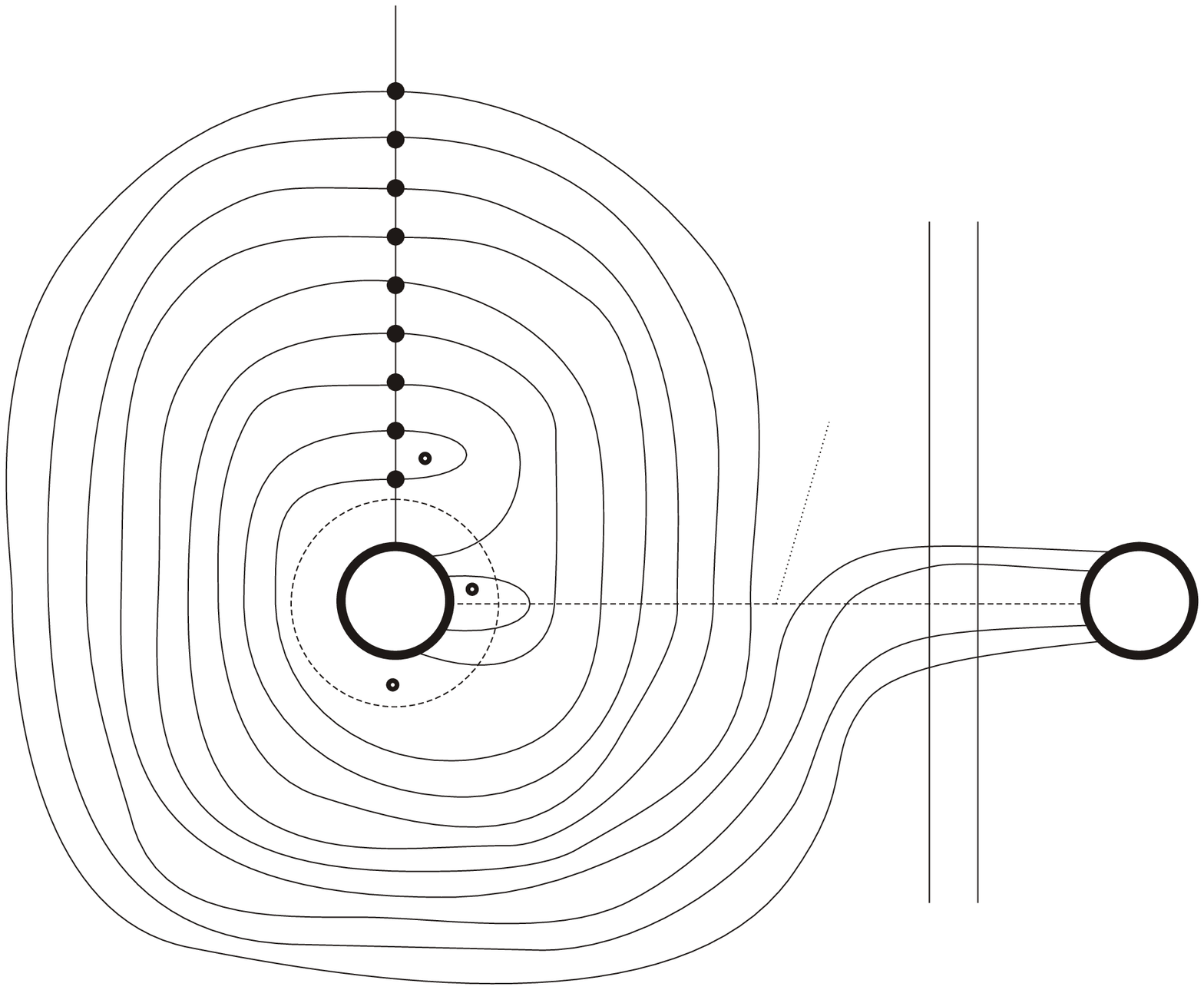}}}

\put(255,201){$\lambda_i$}

\put(250,50){$\gamma_i$}

\put(300,270){$\alpha$}

\put(112,113){$\scriptstyle t_i$}

\put(144,133){$\scriptstyle z_i'$}

\put(133,184){$\scriptstyle w_i$}

\put(110,179){${\scriptstyle x_i^0}$}

\put(110,304){${\scriptstyle x_i^4}$}

\put(112,332){$\alpha_i$}

\put(142,105){$\beta_i$}
\end{picture}
Figure 1\quad The local Heegaard diagram for a $(3,7)$--cable.
\end{center}

The curve $\alpha_i$ has $2(p_i-1)n_i+1$ intersection points with
$\gamma_i$. We define a function
$$S_i\co\alpha_i\cap\gamma_i\to\mathbb Z,$$
which is uniquely characterized up to an overall translation as
follows.

Given $x,y\in\alpha_i\cap\gamma_i$, there are two arcs
$a\subset\alpha_i$ and $b\subset\gamma_i$, both connecting $x$ to
$y$, so that $a-b$ is homologous to a sum of copies of $\alpha$
curves and $\gamma$ curves. Let $\mathcal D$ be such a homology,
then $S_i$ satisfies
$$S_i(x)-S_i(y)=n_{z_i'}(\mathcal D)-n_{w_i}(\mathcal D).$$

We claim that, the values of $S_i$ at the $2(p_i-1)n_i+1$ points
in $\alpha_i\cap\gamma_i$ are mutually different. In order to show
this, we consider a model, which is the $(p,pn+1)$ cable of the
unknot. (In this paragraph, we suppress the subscript $i$ in $p_i$
and $n_i$ for simplicity.) An elementary calculation shows that
the Alexander polynomial of the torus knot $K=T(p,pn+1)$ is
$$\Delta_{K}(T)=T^{-g}(\sum_{j=0}^{(p-1)n}T^{pj}-\sum_{k=0}^{p-2}\sum_{h=0}^{n-1}T^{k(pn+1)+ph+1}),$$
which has exactly $2(p-1)n+1$ terms. Hence the values of $S_i$ at
these $2(p-1)n+1$ points are mutually different. Moreover, if we
denote these points by $x_i^0,x_i^1,\dots,x_i^{2(p-1)n}$, so that
$$S_i(x_i^j)>S_i(x_i^k), \quad\textrm{if}\quad j<k,$$
then $S_i$ satisfies the formula
$$\sum_j(-1)^jT^{S_i(x_i^j)}=T^m\Delta_{K}(T).$$

It is easy to see that the point $x_i^0$ comes from the original
intersection point $\alpha_i\cap\beta_i$, and $x_i^0$ is the
``outermost" point in $\alpha_i\cap\gamma_i$. The point
$x_i^{2n_i}$ is an ``innermost" point in $\alpha_i\cap\gamma_i$,
satisfying
\begin{equation}\label{diffpn}
S_i(x_i^0)-S_i(x_i^{2n_i})=p_in_i.
\end{equation}

Let $A$ be a non-empty subset of $\{1,\dots,l\}$. If the
$\gamma_i$ coordinate of an intersection point $\mathbf
u\in\mathbb T_{\alpha}\cap\mathbb T_{\gamma}$ is supported in a
regular neighborhood of $\lambda_i$ when $i\in A$, and is $x_i^0$
when $i\notin A$, then $\mathbf u$ is called a {\it type--$A$
outermost interior intersection point}. The set of type--$A$
outermost interior intersection points is denoted by
$\mathrm{O}^A$. Given an intersection point $\mathbf x\in\mathbb
T_{\alpha}\cap\mathbb T_{\beta}$, one can associate to it a
corresponding intersection point $\mathbf x'\in\mathbb
T_{\alpha}\cap\mathbb T_{\gamma}$, so that the $\gamma_i$
coordinate of $\mathbf x'$ is $x_i^0$ for all $i=1,\dots,l$. We
then call $\mathbf x'$ an {\it outermost exterior intersection
point}. We can also associate to $\mathbf x$ a corresponding
intersection point $\mathbf x^A\in\mathbb T_{\alpha}\cap\mathbb
T_{\gamma}$, so that the $\gamma_i$ coordinate of $\mathbf x^A$ is
$x_i^{2n_i}$ if $i\in A$, is $x_i^{0}$ if $i\notin A$, and all
other coordinates are the same as the coordinates of $\mathbf x$.

In order to emphasize the dependence of the diagram on $\mathbf
n$, we sometimes put a subscript $(\mathbf n)$ in the notation.
For example, the base points $\mathbf z'$ are denoted by $\mathbf
z'_{(\mathbf n)}$, and the set of type--$A$ outermost interior
intersection points is denoted by $\mathrm{O}_{(\mathbf n)}^A$.

For two different $\mathbf n_1,\mathbf n_2$, there is a natural
1--1 correspondence between $\mathrm{O}^A_{(\mathbf n_1)}$ and
$\mathrm{O}^A_{(\mathbf n_2)}$. Suppose $\mathbf u_{(\mathbf
n_1)}\in\mathrm{O}^A_{(\mathbf n_1)}$, then the corresponding
point in $\mathrm{O}^A_{(\mathbf n_2)}$ is denoted by $\mathbf
u_{(\mathbf n_2)}$.

Let $$\begin{array}{ll} \mathfrak h_{\mathbf w,\mathbf
z}\co&\mathbb T_{\alpha}\cap\mathbb
T_{\beta}\to\homo^2(Y,L;\mathbb Q),\\
\mathfrak h_{\mathbf w,\mathbf t}\co&\mathbb T_{\alpha}\cap\mathbb
T_{\gamma}\to\homo^2(Y,L;\mathbb Q),\\
\mathfrak h_{\mathbf w,\mathbf z'}\co&\mathbb
T_{\alpha}\cap\mathbb T_{\gamma}\to\homo^2(Y,C(L);\mathbb Q),
\end{array}$$ be the affine maps defined in Subsection
\ref{sectfil}.

The following observation is important:

\begin{lem}\label{filconst}
Fix a point $\mathbf x\in\mathbb T_{\alpha}\cap\mathbb T_{\beta}$,
a set $A\subset\{1,\dots,l\}$ and a point $\mathbf u_{(\mathbf n
)}\in \mathrm{O}^A_{(\mathbf n)}$. Then
$$\mathfrak h_{\mathbf w,\mathbf z'_{(\mathbf n)}}(\mathbf x^A_{(\mathbf n)})-\mathfrak h_{\mathbf w,\mathbf z'_{(\mathbf n)}}(\mathbf u_{(\mathbf n)})$$
is a constant independent of $\mathbf n$.
\end{lem}
\begin{proof}
Given two $l$--tuples $\mathbf n_1,\mathbf n_2$. Without loss of
generality, we can assume $\mathbf n_1<\mathbf n_2$, that is,
every coordinate of $\mathbf n_1$ is less than or equal to the
corresponding coordinate of $\mathbf n_2$, and at least one
equality does not hold. Suppose $\mathcal D_{(\mathbf n_1)}$ is
the domain which gives a homology between $k\omega(\mathbf
x^A_{(\mathbf n_1)},\mathbf u_{(\mathbf n_1)})$ and a sum of some
copies of $\alpha$ curves and $\gamma_{(\mathbf n_1)}$ curves.
Then  we can get a domain $\mathcal D_{(\mathbf n_2)}$ by
performing finger moves to $\mathcal D_{(\mathbf n_1)}$, so that
$\mathcal D_{(\mathbf n_2)}$ is the corresponding domain for
$\mathbf x^A_{(\mathbf n_2)},\mathbf u_{(\mathbf n_2)}$.

When $i\in A$, there is an arc $\zeta$, supported in the $i$--th
spiral, connecting $w_i$ to $z'_{i(\mathbf n)}$, and
$n_{z'_{i(\mathbf n)}}(\mathcal D_{(\mathbf n)})-n_{w_i}(\mathcal
D_{(\mathbf n)})$ is calculated by the algebraic intersection
number of $\zeta$ with $\partial\mathcal D_{(\mathbf n)}$. Since
the $\gamma_i$ coordinate of $\mathbf u_{(\mathbf n)}$ is
supported in $Nd(\lambda_i)$, and the $\gamma_i$ coordinate of
$\mathbf x^A_{(\mathbf n)}$ is an ``innermost" point $x_{(\mathbf
n)}^{2n_i}$, it is easy to see that the finger moves do not change
the algebraic intersection number of $\zeta$ with
$\partial\mathcal D_{(\mathbf n)}$.

When $i\notin A$, the $\gamma_i$ coordinates of $\mathbf
x^A_{(\mathbf n)}$ and $\mathbf u_{(\mathbf n)}$ are both
$x_{i(\mathbf n)}^0$. It is easy to see that the finger moves do
not change $n_{z'_{i(\mathbf n)}}(\mathcal D_{(\mathbf
n)})-n_{w_i}(\mathcal D_{(\mathbf n)})$. Thus our desired result
holds by Lemma \ref{frakH}.
\end{proof}

\begin{lem}\label{h=h=h}
Given $\mathbf x,\mathbf y\in\mathbb T_{\alpha}\cap\mathbb
T_{\beta}$, we have
$$
\mathfrak h_{\mathbf w,\mathbf z}(\mathbf x)-\mathfrak h_{\mathbf
w,\mathbf z}(\mathbf y)=\mathfrak h_{\mathbf w,\mathbf t}(\mathbf
x')-\mathfrak h_{\mathbf w,\mathbf t}(\mathbf y')=j^*(\mathfrak
h_{\mathbf w,\mathbf z'}(\mathbf x')-\mathfrak h_{\mathbf
w,\mathbf z'}(\mathbf y')).
$$
\end{lem}
\begin{proof}
It is obvious that
$$\mathfrak h_{\mathbf w,\mathbf z}(\mathbf x)-\mathfrak h_{\mathbf
w,\mathbf z}(\mathbf y)=\mathfrak h_{\mathbf w,\mathbf t}(\mathbf
x')-\mathfrak h_{\mathbf w,\mathbf t}(\mathbf y').$$

Suppose $\mathcal D$ is a domain, $\partial\mathcal D$ is the sum
of $k\omega(\mathbf x,\mathbf y)$ and some copies of $\alpha$ and
$\beta$ curves. Then after applying $\mathbf p$--fold finger moves
to $\mathcal D$, we get a domain $\mathcal D'$, so that
$\partial\mathcal D'$ is the sum of $k\omega(\mathbf x',\mathbf
y')$ and some copies of $\alpha$ and $\gamma$ curves. It is not
hard to show
$$n_{z_i'}(\mathcal D')-n_{w_i}(\mathcal D')=p_i(n_{t_i}(\mathcal D')-n_{w_i}(\mathcal D')).$$
Hence $$j^*(\mathfrak h_{\mathbf w,\mathbf z'}(\mathbf
x')-\mathfrak h_{\mathbf w,\mathbf z'}(\mathbf y'))=\mathfrak
h_{\mathbf w,\mathbf t}(\mathbf x')-\mathfrak h_{\mathbf w,\mathbf
t}(\mathbf y')$$ by Lemma \ref{frakH} and the fact that
$j^*(p_i\mathrm{PD}([\mu_i']))=\mathrm{PD}([\mu_i]).$
\end{proof}

Suppose $\mathbf x^j,\mathbf x^k\in\mathbb T_{\alpha}\cap\mathbb
T_{\gamma}$ differ only at the $\gamma_i$ coordinate, where the
coordinate of $\mathbf x^j$ is $x^j_i$, and the coordinate of
$\mathbf x^k$ is $x^k_i$. From the definition of $S_i$ and Lemma
\ref{frakH}, we conclude that
\begin{equation}\label{diffjk}
\mathfrak h_{\mathbf w,\mathbf z'}(\mathbf x^j)-\mathfrak
h_{\mathbf w,\mathbf z'}(\mathbf
x^k)=(S_i(x_i^j)-S_i(x_i^k))\mathrm{PD}([\mu_i']).
\end{equation}

Now we fix an integral class $h\in\homo_2(Y,L;\mathbb Q)$, which
satisfies $h\cdot[\mu_i]>0$ for each $i$. Suppose $F\subset
M=Y-\mathrm{int}(Nd(L))$ is a surface representing $h$, $F$ has no
sphere components, and $\chi(F)$ is maximal among all such
surfaces. We can assume $\partial F\cap\partial Nd(K_i)$ consists
of parallel oriented circles. Then $\partial F\cap\partial
Nd(K_i)$ is a torus link $T(P_i,Q_i)$, with respect to the frame
specified by $\Sigma$.

From now on, we assume $p_i/P_i$ is an integer independent of $i$,
say, $p_i=mP_i$. Then $C_{\mathbf p,\mathbf q}(L)$ is a
null-homologous link. In fact, a minimal genus Seifert surface
$F'$ for $C(L)$ can be obtained as follows. Inside the cable space
$Nd(K_i)-\mathrm{int}(Nd(C_{p_i,q_i}(K_i)))$, one can choose a
properly embedded, Thurston norm minimizing surface $G_i$, so that
$\partial G_i\cap\partial Nd(K_i)$ is the torus link
$T(mP_i,mQ_i)$, and $\partial G_i\cap\partial
Nd(C_{p_i,q_i}(K_i))$ is a longitude of $C_{p_i,q_i}(K_i)$. Then
$F'$ is the union of $G_1,\dots,G_l$ and $m$ parallel copies of
$F$. A standard argument in 3--dimensional topology shows that
$F'$ is a minimal genus Seifert surface for $C(L)$. Let
$h'=[F']\in\homo_2(Y,C(L))$.

Recall the function
$$\mathcal F^h_{\mathbf w,\mathbf z}\co\mathbb T_{\alpha}\cap\mathbb T_{\beta}\to\mathbb Q$$
is defined as
$$\mathcal F^h_{\mathbf w,\mathbf z}(\mathbf x)=\langle\mathfrak h_{\mathbf w,\mathbf z}(\mathbf x),h\rangle.$$
Then $\mathcal F^h_{\mathbf w,\mathbf z}$ specifies an Alexander
$\mathbb Q$--grading on $\widehat{CFL}(Y,L)$. We also equip
$\widehat{CFL}(Y,C(L))$ with the $\mathbb Q$--grading defined by
$\mathcal F^{h'}_{\mathbf w,\mathbf z'}$. Let
$\widehat{HFL}(Y,L,\textrm{topmost})$,
$\widehat{HFL}(Y,C(L),\textrm{topmost})$ be the topmost nontrivial
terms in $\widehat{HFL}(Y,L)$ and $\widehat{HFL}(Y,C(L))$, with
respect to these $\mathbb Q$--gradings, respectively. If the
grading of $\mathbf x$ is no more than the grading of $\mathbf y$,
then we denote as $\mathbf x\preceq\mathbf y$.

Given $i\in\{1,\dots,l\}$, suppose $\mathbf x^j,\mathbf
x^k\in\mathbb T_{\alpha}\cap\mathbb T_{\gamma}$ are two points
differing only at the $\gamma_i$ component, where their components
are $x_i^j,x_i^k$, respectively. By (\ref{diffjk}), we have
\begin{equation}\label{diffjkF}
\mathcal F^{h'}_{\mathbf w,\mathbf z'}(\mathbf x_i^j)-\mathcal
F^{h'}_{\mathbf w,\mathbf z'}(\mathbf
x_i^k)=S_i(x_i^j)-S_i(x_i^k).
\end{equation}
Moreover, if $\mathbf x,\mathbf y\in\mathbb T_{\alpha}\cap\mathbb
T_{\beta}$ are two intersection points, then by Lemma \ref{h=h=h}
and the construction of $F'$, we have
\begin{equation}\label{diff''F}
\mathcal F^{h'}_{\mathbf w,\mathbf z'}(\mathbf x')-\mathcal
F^{h'}_{\mathbf w,\mathbf z'}(\mathbf y')=m(\mathcal
F^{h}_{\mathbf w,\mathbf z}(\mathbf x)-\mathcal F^{h}_{\mathbf
w,\mathbf z}(\mathbf y)).
\end{equation}

\begin{prop}\label{cabletop}
With the notation as above, when the winding number $\mathbf n$ is
sufficiently large, the following equality holds
$$\widehat{HFL}(Y,L,\text{\rm topmost})\cong\widehat{HFL}(Y,C(L),\text{\rm topmost}).$$
Moreover, suppose $\mathbf x$ is one of the generators of
$\widehat{CFL}(Y,L,\text{\rm topmost})$, then $\mathbf x'$ is one
of the generators of $\widehat{CFL}(Y,C(L),\text{\rm topmost})$.
\end{prop}
\begin{proof}
Let $\{\mathbf x_1\dots\mathbf x_r\}\subset\mathbb
T_{\alpha}\cap\mathbb T_{\beta}$ be the generating set of
$\widehat{CFL}(Y,L,\textrm{topmost})$. If $\mathbf v$ is not an
exterior intersection point, we change all of its $\gamma_i$
components which are supported in $Nd(\beta_i)$ to $x^0_i$, so as
to get a type--$A$ outermost interior intersection point
$\widetilde{\mathbf v}$, for some $A$. By (\ref{diffjkF}) we have
\begin{equation}\label{v<v}
\mathbf v\preceq\widetilde{\mathbf v}
\end{equation}

Let $C_1$ be a lower bound of $\mathcal F^{h'}_{\mathbf w,\mathbf
z'}(\mathbf x_1^A)-\mathcal F^{h'}_{\mathbf w,\mathbf z'}(\mathbf
u)$ for all nonempty $A\subset\{1,\dots,l\}$, and all type--$A$
outermost interior intersection points $\mathbf u$.
Lemma~\ref{filconst} enables us to choose $C_1$ to be a constant
independent of $\mathbf n$.

Let $\mathbf n$ be sufficiently large so that
\begin{equation}\label{choiceN}
p_in_i+C_1>0, \qquad\text{for any}\quad i=1,\dots,l.
\end{equation}

Let $\widehat{CFL}_{\preceq}$ be the summand of
$\widehat{CFL}(Y,C(L))$, which consists of all the elements with
grading no lower than the grading of $\mathbf
 x'_1$. By (\ref{diffpn}),
(\ref{diffjkF}), (\ref{v<v}) and (\ref{choiceN}), we find that the
generators of $\widehat{CFL}_{\preceq}$ are all exterior
intersection points. The differential on $\widehat{CFL}_{\preceq}$
counts holomorphic disks away from $\mathbf w,\mathbf z'$, denoted
by $\partial_{\mathbf w,\mathbf z'}$.

The base points $\mathbf t$ give an extra filtration to
$\widehat{CFL}_{\preceq}$. It is easy to see that if a holomorphic
disk $\phi$ connects two exterior points $\overline{\mathbf y}_1$
to $\overline{\mathbf y}_2$, and $\phi$ avoids $\mathbf w,\mathbf
z',\mathbf t$, then the $\gamma_i$ components of
$\overline{\mathbf y}_1$ and $\overline{\mathbf y}_2$ coincide for
all $i$. Thus $\phi$ corresponds to a holomorphic disk connecting
$\mathbf y_1$ to $\mathbf y_2$, which avoids $\mathbf w,\mathbf
z$. Here $\mathbf y_j\in\mathbb T_{\alpha}\cap\mathbb T_{\beta}$,
($j=1,2$,) is an intersection point whose components coincide with
the components of $\overline{\mathbf y}_j$, except the $\beta_i$
components.

Hence the chain complex
$(\widehat{CFL}_{\preceq},\partial_{\mathbf w,\mathbf z',\mathbf
t})$ is the direct sum of summands in the form of
$\widehat{CFL}_{\mathbf j,d}$, where here $\widehat{CFL}_{\mathbf
j,d}$ is generated by the exterior intersection points
$\overline{\mathbf y}\in\mathbb T_{\alpha}\cap\mathbb T_{\gamma}$,
which satisfy that the $\gamma_i$ component of $\mathbf y$ is
$x_i^{j_i}$, and the grading difference between $\overline{\mathbf
y}$ and $\mathbf x'_1$ is $d\ge0$. Moreover, by (\ref{diffjkF})
and (\ref{diff''F}), the homology of $(\widehat{CFL}_{\mathbf
j,d},\partial_{\mathbf w,\mathbf z',\mathbf t})$ is isomorphic to
the homology of some summand of $\widehat{CFL}(Y,L)$, at some
grading no less than the grading of $\mathbf x_1$.

Since $\mathbf x_1$ lies in the topmost nontrivial Alexander
$\mathbb Q$--grading of $\widehat{HFL}(Y,L)$, we find that
$\widehat{HFL}_{\mathbf j,d}$ is nontrivial if and only if
$(\mathbf j,d)=(\mathbf 0,0)$, and
$$\widehat{HFL}_{\mathbf 0,0}\cong\widehat{HFL}(Y,L,\text{\rm topmost}).$$ There is a
spectral sequence which starts from
$(\widehat{CFL}_{\preceq},\partial_{\mathbf w,\mathbf z',\mathbf
t})$, and converges to
$\homo(\widehat{CFL}_{\preceq},\partial_{\mathbf w,\mathbf z'})$.
Since the $E^2$ term is only supported in one filtration level, we
must have
$$\homo(\widehat{CFL}_{\preceq},\partial_{\mathbf w,\mathbf z'})\cong\widehat{HFL}_{\mathbf 0,0}.$$
Thus
$$\widehat{HFL}(Y,C(L),\text{\rm topmost})\cong\widehat{HFL}(Y,L,\text{\rm topmost}).$$

The last statement of this proposition is obvious from the proof.
\end{proof}

Our next task is to determine the absolute position of the topmost
grading in $\widehat{HFL}(Y,L)$. For this purpose, we will define
two functions
$$\mathcal F_1,\mathcal F_2\co\mathbb T_{\alpha}\cap\mathbb T_{\beta}\to\mathbb Q$$
as follows.

The cable space $Nd(K_i)-\mathrm{int}(Nd(C_{p_i,q_i}(K_i)))$
fibers over the circle, with fiber $G_i$. Let $u_i$ be a vector
field on the cable space, which is transverse to the fibers
everywhere, and the orientation of $u_i$ is {\bf opposite} to the
orientation induced by the orientation of the fibers. Moreover,
let the restriction of $u_i$ on the boundary tori be the canonical
translation invariant vector field.

Given $\mathbf x\in\mathbb T_{\alpha}\cap\mathbb T_{\beta}$, let
$\mathfrak r_1=\underline{\mathfrak s}_{\mathbf w,\mathbf
z'}(\mathbf x')\in\relspin(Y,C(L))$. The function $\mathcal F_1$
is defined as follows:
$$\mathcal F_1(\mathbf x)=\langle \mathfrak H_{C(L)}(\mathfrak r_1),h'\rangle,$$
where here $\mathfrak H_{C(L)}$ is the affine map defined in
Subsection \ref{sectfil}, for the pair $(Y,C(L))$.

Note that $\underline{\mathfrak s}_{\mathbf w,\mathbf z}(\mathbf
x)=\underline{\mathfrak s}_{\mathbf w,\mathbf t}(\mathbf x')$ is a
relative Spin$^c$ structure on $Y-\mathrm{int}(Nd(L))$, we can
extend it to a relative Spin$^c$ structure on
$Y-\mathrm{int}(Nd(C(L)))$ by the vector fields $u_1,\dots,u_l$.
We denote this new Spin$^c$ structure by $\mathfrak r_2$. Now let
$$\mathcal F_2(\mathbf x)=\langle \frac{c_1(\mathfrak r_2)-\sum_{i=1}^l\mathrm{PD}([\mu_i])}2,h'\rangle.$$

In summary, $\mathcal F_1$ and $\mathcal F_2$ can be factorized as
follows:
\begin{eqnarray*}
\mathcal F_1\co&&\mathbb T_{\alpha}\cap\mathbb T_{\beta}\to\mathbb
T_{\alpha}\cap\mathbb T_{\gamma}\quad\;\;\to\relspin(Y,C(L))\to\mathbb Q,\\
\mathcal F_2\co&&\mathbb T_{\alpha}\cap\mathbb
T_{\beta}\to\relspin(Y,L)\to\relspin(Y,C(L))\to\mathbb Q.
\end{eqnarray*}

Let $\max\mathcal F_j$ be the maximal value of $\mathcal
F_j(\mathbf x)$, where $\mathbf x$ runs over the nontrivial
filtration levels of $\widehat{HFL}(Y,L)$.

From Lemma \ref{h=h=h}, we can conclude that
\begin{equation}\label{F1F2}
\mathcal F_1(\mathbf x)-\mathcal F_1(\mathbf y)=\mathcal F
_2(\mathbf x)-\mathcal F_2(\mathbf y).
\end{equation}

In fact, we can prove the stronger
\begin{prop}
Given $\mathbf x\in\mathbb T_{\alpha}\cap\mathbb T_{\beta}$, then
the following equality holds:
\begin{equation}\label{F1F2eqn}
\mathcal F_2(\mathbf x)=\mathcal F_1(\mathbf
x)+\sum_{i=1}^l\frac{p_i-1}2.\end{equation}
\end{prop}
\begin{proof}We could prove (\ref{F1F2eqn}) by examining the
relative Spin$^c$ structures carefully. But we would rather argue
via a model computation.

Let $f$ be a Morse function corresponding to the Heegaard diagram
$(\Sigma,\mbox{\boldmath${\alpha}$},\mbox{\boldmath${\gamma}$})$.
Following the definitions of $\mathfrak r_1,\mathfrak r_2$ and the
construction in Subsection \ref{sectdg}, we can construct two
vector fields $v_1,v_2$ on $Y$, representing the two relative
Spin$^c$ structures.

We note that $v_1,v_2$ are equal outside a regular neighborhood of
the flowlines $\gamma_{\mathbf w},\gamma_{\mathbf
z'},\gamma_{\mathbf t}$. And the difference of $v_1$ and $v_2$
inside $Nd(\gamma_{\mathbf w}\cup\gamma_{\mathbf
z'}\cup\gamma_{\mathbf t})$ depends only on the $2l$ torus link
types $(p_i,q_i),(mP_i,mQ_i)$, $i=1,\dots,l$. Moreover, one can
isotope $F'$ so that $F'\cap Nd(\gamma_{\mathbf
w}\cup\gamma_{\mathbf z'}\cup\gamma_{\mathbf t}))$ depends only on
$(\mathbf p,\mathbf q)$ and $(m\mathbf P,m\mathbf Q)$. So we only
need to verify (\ref{F1F2}) for some models.

Let $d_i=\gcd(P_i,Q_i)$, $P_i=d_i P_i',Q_i=d_i Q_i'$. Consider the
knot $O_{P_i'/Q_i'}$ in $L(P_i',Q_i')$. There is an essential disk
$D$ properly embedded in the complement of $K=O_{P_i'/Q_i'}$.
$\partial D$ is the torus knot $T(P_i',Q_i')$ in $\partial Nd(K)$.
Let $F_0$ be the union of $d_i$ copies of $D$.

Let $C(K)$ be the $(p_i,q_i)$--cable of $K$, $G_0\subset
Nd(K)-\mathrm{int}(Nd(C(K)))$ is a surface, such that $\partial
G_0$ consists of the torus link $T(mP_i,mQ_i)$ and a longitude of
$C(K)$, and $Nd(K)-\mathrm{int}(Nd(C(K)))$ fibers over the circle
with fiber $G_0$. $F_0'$ is the union of $G_0$ and $m$ parallel
copies of $F_0$. $F_0'\cap Nd(\gamma_{\mathbf
w}\cup\gamma_{\mathbf z'}\cup\gamma_{\mathbf t})$ is in the same
pattern as before.

From Lemma \ref{LensKnot}, we know that for the pair
$(L(P_i',Q_i'),O_{P_i'/Q_i'})$
\begin{eqnarray*}
\max\mathcal F_2&=&\frac12(md_i(2P_i'-1)-\chi(G_0)-p_i)\\
&=&\frac12(p_i-md_i-\chi(G_0)).
\end{eqnarray*}

The knot $C(K)$ is null-homologous, so we can apply Proposition
\ref{LinkNorm1} and Proposition \ref{cabletop} to show that
\begin{eqnarray*}\max\mathcal F_1&=&\frac{1-\chi(F_0')}2\\
&=&\frac12(1-md_i-\chi(G_0)).
\end{eqnarray*}

So we get
$$\mathcal F_2(\mathbf x)-\mathcal F_1(\mathbf x)=\frac{p_i-1}2$$
in this case. Hence (\ref{F1F2eqn}) holds in general.
\end{proof}

\begin{proof}[Proof of Theorem \ref{LinkNorm}]We first prove the case where
$\partial M$ is incompressible, thus we only need to prove the
statement in Remark \ref{Thurston}. By the continuity and
linearity of Thurston norm, it suffices to prove the theorem for
the integral classes $h\in\homo_2(Y,L)$ which satisfy
$$h\cdot[\mu_i]>0$$
for all $i$. With the notation as before, consider the $(\mathbf
p,\mathbf q)$--cable $C(L)$ of $L$. Here we choose $p_i=P_i$, and
$q_i=p_in_i+1$. Let $\mathbf n=(n_1,\dots,n_l)$ be sufficiently
large. Since $\partial M$ is incompressible, there exists a
surface $F$ representing $h$, and $|\chi(F)|=x(h)$. Construct the
surfaces $G_i$,$F'$ as before.

By (\ref{F1F2eqn}), Proposition \ref{cabletop} and Proposition
\ref{LinkNorm1},
\begin{eqnarray*}
&&\max\{\langle \mathfrak H(\mathfrak
r),h\rangle|\:\widehat{HFL}(Y,L,\mathfrak
r)\ne0\}\\
&=&\max\mathcal
F_2+\frac12\sum_{i=1}^l\chi(G_i)\\
&=&\max\mathcal
F_1+\sum_{i=1}^l\frac{p_i-1}2+\frac12\sum_{i=1}^l\chi(G_i)\\
&=&\frac12(l-\chi(F'))+\frac12\sum_{i=1}^l(p_i-1)+\frac12\sum_{i=1}^l\chi(G_i)\\
&=&\frac12(\sum_{i=1}^l|[F]\cdot[\mu_i]|-\chi(F)).
\end{eqnarray*}
This finishes the proof in the case when $\partial M$ is
incompressible.

If $\partial M$ is compressible, say, $\partial Nd(K_1)$ is
compressible. We can compress this boundary torus to get a
separating sphere, which splits off a lens space summand from $Y$,
and $K_1$ is a knot $O_{p/q}$ in this summand. Let $L'=L-K_1$,
$Y=Y'\#L(p,q)$. It is easy to see that, if $h\in\homo_2(Y,L)$ is
an integral class, and $F\subset Y-\mathrm{int}(Nd(L))$ realizes
$\chi(h)$, then $F$ is the disjoint union of some disks in
$L(p,q)-\mathrm{int}(Nd(O_{p/q}))$ and a surface $F'\subset
Y'-\mathrm{int}(Nd(L'))$. We can make use of
Proposition~\ref{splitlink} and Lemma~\ref{LensKnot} to reduce our
problem to $L'$. Now the proof of our theorem can be finished by
induction on $|L|$.
\end{proof}

\end{document}